\documentclass[11pt]{article}

\usepackage{amsmath,amssymb,amsthm}
\usepackage{mathrsfs}
\usepackage{fullpage}
\usepackage{hyperref}
\usepackage[alphabetic]{amsrefs}

\newcommand{\ceiling}[1]{\ensuremath{\lceil #1 \rceil}}

\newcommand{\RR}{\ensuremath{{\mathbb R}}}
\newcommand{\ZZ}{\ensuremath{{\mathbb Z}}}

\newcommand{\TT}{\ensuremath{{\mathbb T}}}
\newcommand{\sN}{\mathscr{N}}

\renewcommand{\vec}{\overline}
\newcommand{\vt}{\vec\omega}
\newcommand{\va}{\vec\alpha}

\newcommand{\Expect}[2]{\ensuremath{{\mathbb E}_{#1}\left[#2\right]}}
\newcommand{\Prob}[2]{\ensuremath{{\mathbb P}_{#1}\left[#2\right]}}

\DeclareMathOperator{\vol}{{\bf vol}}

\DeclareMathOperator{\BOX}{{\textsc{Box}}}
\DeclareMathOperator{\annuli}{{\textsc{Annuli}}}
\DeclareMathOperator{\ball}{{\textsc{Ball}}}

\DeclareMathOperator{\diam}{{diam}}
\DeclareMathOperator{\type}{{\textsc{Type}}}

\newtheorem{theorem}{Theorem}

\newtheorem{lemma}{Lemma}

\newtheorem{corollary}{Corollary}

\title{Thick subsets that do not contain arithmetic progressions}
\author{Kevin O'Bryant \\ City University of New York, College of Staten Island and The Graduate Center}
\date{\today}

\begin{document}

\maketitle
\begin{abstract}
We adapt the construction of subsets of $\{1,2,\dots,N\}$ that contain no $k$-term arithmetic progressions to give a relatively thick subset of an arbitrary set of $N$ integers. Particular examples include a thick subset of $\{1,4,9,\dots,N^2\}$ that does not contain a 3-term AP, and a positive relative density subset of a random set (contained in $\{1,2,\dots,n\}$ and having density $c n^{-1/(k-1)}$) that is free of $k$-term APs.
\end{abstract}

\section{Introduction}

For a  finite set $\sN$ (whose cardinality we denote by $N$), set $r_k(\sN)$ to be the largest possible size of a subset of $\sN$ that contains no $k$-term arithmetic progressions ($k$-APs). A little-known result of Koml\'{o}s, Sulyok, and Szemer\'{e}di~\cite{KSSz} implies that
    \begin{equation}\label{KSSz}
    r_k(\sN) \geq C \,r_k([N]),
    \end{equation}
for an explicit positive constant $C$. Abbott~\cite{Abbott} reports that their proof gives $C=2^{-15}$, and indicates some refinements that yield $C=1/34$.

In this work, we focus on the situation when $\sN$ itself has few solutions and we can give bounds on $r_k(\sN)$ that are much stronger than those implied by Eq.~\eqref{KSSz} and the currently best bounds on $r_k([N])$. In particular, we adapt the Behrend-type construction~\cite{OBryant} of subsets of $[N]$ without $k$-APs to arbitrary finite sets $\sN$. We draw particular attention to subsets of the squares and to subsets of random sets. As the statement of our theorem requires some notation and terminology, we first give two corollaries.

Our first corollary brings attention to the fact that while the squares contain many 3-APs, they also contain unusually large subsets that do not. Here and throughout this paper, $\exp(x)=2^x$ and $\log x=\log_2(x)$. For comparison, $r_3([N]) \geq N\exp(-2\sqrt{2}\sqrt{\log N}+\frac14\log\log N)$.
\begin{corollary}\label{cor:squares}
There is an absolute constant $C>0$ such that for every $N$ there is a subset of $\{1,4,9,\dots,N^2\}$ with cardinality at least
\[C \,\cdot\,N\,\cdot\, \exp\big({-2\sqrt{2}\sqrt{\log\log N} + \frac14 \log\log\log N}\big)\]
that does not contain any 3-term arithmetic progressions.
\end{corollary}

Our second corollary identifies sets that have subsets with no $k$-APs and with {\em positive} relative density .
\begin{corollary}\label{cor:random}
For every real $\psi$ and integer $k\geq 3$, there is a real $\delta>0$ such that every sufficiently large $\sN\subseteq \ZZ$ that has fewer than $\psi|\sN|$ arithmetic progressions of length $k$ contains a subset that is free of $k$-term arithmetic progressions and has relative density at least $\delta$. In particular, for each $\delta>0$, if $n$ is sufficiently large and $\sN\subseteq \{1,2,\dots,n\}$ is formed by including each $k$ independently with probability $c n^{-1/(k-1)}>0$, then with high probability $\sN$ contains a subset $A$ with relative density $\delta$ and no $k$-term arithmetic progressions.
\end{corollary}

The structure of the proof requires us to consider a generalization of arithmetic progressions. A $k$-term $D$-progression is a nonconstant sequence $a_1,\dots,a_k$ whose $(D+1)$-st differences are all zero:
    \[\sum_{i=0}^{D+1} (-1)^i \binom{D+1}{i} a_{i+v} = 0,\qquad (1\leq v \leq k-D-1).\]
Equivalently, $a_1,\dots,a_k$ is a $k$-term $D$-progression if there is a nonconstant polynomial $Q(j)$ with degree at most $D$ and $Q(i)=a_i$ for $i\in[k]$.  Clarifying examples of $5$-term $2$-progressions of integers are $1,2,3,4,5$ (from $Q(j)=j$), and $4,1,0,1,4$ (from $Q(j)=(j-3)^2$), and $1,3,6,10,15$ (from $Q(j)=\frac12 j+\frac12 j^2$). Let $Q(j)=\sum_{i=0}^{D'} q_i j^i$ be a polynomial with degree $D'\geq1$, so that $Q(1),Q(2),\dots,Q(k)$ is a $k$-term $D$-progression for all $D\geq D'$. The quantity $D'!q_{D'}$, which is necessarily nonzero, is called the difference of the sequence, and $(D',Q(1),D'!q_{D'})$ is the type of the sequence. Note that different progressions can have the same type: both $1,4,9,16,25$ and $1,5,11,19,29$ have type $(2,1,2)$. For any set $\sN$, we let $\type_{k,D}(\sN)$ be the number of types of $k$-term $D$-progressions contained in $\sN$. The proof of \cite{OBryant}*{Lemma 4} shows that $\type_{k,D}(\sN) \ll |\sN| \diam(\sN)$. Since the type of a $k$-term $D$-progression is determined by its first $D+1$ elements, we also have $\type_{k,D}(\sN) \leq {N}^{D+1}$. We define
    \[r_{k,D}(\sN):= \max_{A\subseteq \sN} \left\{ |A| \colon A\text{ does not contain any $k$-term $D$-progressions}\right\}\]
and recall the lower bound proved in~\cite{OBryant}:
    \begin{equation}\label{equ:OBryant}
    \frac{r_{k,D}([N])}{N} \geq C \exp\left(-n2^{(n-1)/2} D^{(n-1)/n} \sqrt[n]{\log N}+ \frac1{2n}\log\log N\right).
    \end{equation}

We can now state our main theorem.
\begin{theorem}\label{thm:main}
Let $k\geq 3, n\ge 2, D\ge 1$ be integers satisfying $k>2^{n-1}D$.  Let $\Psi(N)$ be any function that is at least 2. There is a constant $C=C(k,D,\Psi)$ such that for all $\sN\subseteq\ZZ$ with $\type_{k,D}(\sN) \leq N\Psi(N)$ (where $N:=|\sN|$)
    \[
    \frac{r_{k,D}(\sN)}{N}  \geq C \exp\left(-n 2^{(n-1)/2} D^{(n-1)/n} \sqrt[n]{\log \Psi(N)}+ \frac{1}{2n}\log\log \Psi(N)\right).
    \]
\end{theorem}

Corollary~\ref{cor:random} is now straightforward: set $D=1$ and $\Psi(N)=\max\{\psi,2\}$ and take
    \[
    \delta = \exp\left(-n 2^{(n-1)/2} \sqrt[n]{\log C}+ \frac{1}{2n}\log\log C\right),
    \]
to arrive at the first sentence. Considering the random set $\sN$ described in the second sentence of Corollary~\ref{cor:random}, for each pair $(a,a+d)$ of elements of $\sN$ the likelihood of the next $k-2$ elements $a+2d,\dots,a+(k-1)d$ of the arithmetic progression being in $\sN$ is $(c n^{-1/(k-1)})^{k-2}$. Consequently, the expected number of $k$-term arithmetic progressions in $\sN$ is
    \[
    \binom{n}{2} (n^{-1/(k-1)})^{k-2} \leq \frac{c^{k-2}}2 n^{k/(k-1)},
    \]
and the expected size of $\sN$ is $N=n \cdot c n^{-1/(k-1)}=c n^{k/(k-1)}$. We can take $\Psi(N)$ to be a constant with high probability, and so Corollary~\ref{cor:random} follows from Theorem~\ref{thm:main}.

Corollary~\ref{cor:squares} is only a bit more involved. It is known (perhaps since Fermat, see \cites{Conrad1,Conrad2,Alf,Brown,Fogarty,KhanKwong,2000clicks} for a history and for the results we use here) that while the squares do not contain any 4-term arithmetic progressions, the 3-term arithmetic progressions $a^2,b^2,c^2$ are parameterized by
    \begin{equation*}
    a = u(2 s t-s^2+t^2),
    b = u(s^2+t^2),
    c = u(2 s t+s^2-t^2),
    \end{equation*}
with $s,t,u\geq1$ and $\gcd(s,t)=1$. Merely observing that $s,t,u\geq 1, b\leq N$ yields that there are fewer than $2\pi N\log N$ triples $(s,t,u)$ with $a,b,c$ in $[N]$, i.e.,
    \[
    \type_{3,1}(\{1,4,9,\dots,N^2\}) \leq 2\pi N \log N.
    \]
Now, setting $k=1,n=2,D=1,\Psi(N)=2\pi\log N$ in Theorem~\ref{thm:main} produces Corollary~\ref{cor:squares}.

Section~\ref{sec:overview} gives a short outline of the construction behind Theorem~\ref{thm:main}, which is given in greater detail in Section~\ref{sec:Proof}. We conclude in Section~\ref{sec:questions} with some unresolved questions.

\section{Overview of construction proving Theorem~\ref{thm:main}}\label{sec:overview}
Throughout this work we fix three integers, $k\geq 3$, $n\geq 2$, $D\geq 1$, that satisfy $k>2^{n-1}D$; in other words, one may take $n=\ceiling{\log(k/D)}$.

In this section, we outline the construction, suppressing as much technical detail as possible. In the following sections, all definitions are made precisely and all arguments are given full rigor.

Fix $\Psi(N)$, and take $\sN\subseteq\ZZ$ with $|\sN|=N$, and so that $\sN$ contains less than $N\Psi(N)$ types of $k$-term $D$-progressions. The parameters $N_0,d,\delta$ are chosen at the end for optimal effect.

Let $A_0=R_{k,2D}(N_0)$ be a subset of $[N_0]$ without $k$-term $2D$-progressions, and
    \[|A_0| = r_{k,2D}(N_0).\]
Consider $\vt,\va$ in $\TT^d$ (we average over all choices of $\vt,\va$ later in the argument), and set
    \[A := \{ a \in \sN \colon a\vt+\va \bmod \vec 1 = \langle x_1,\dots,x_d\rangle, |x_i|< 2^{-D-1}, \sum x_i^2 \in \annuli\},\]
where $\annuli$ is a union of thin annuli in $\RR^d$ with thickness $\delta$ whose radii are affinely related to elements of $A_0$.
Set
    \[T:=\{a \in A \colon \text{there is a $k$-term $D$-progression in $A$ starting at $a$ } \}.\]
Then $A\setminus T$ is free of $k$-term $D$-progressions, and so $r_{k,D}(\sN)\geq |A\setminus T|=|A|-|T|$, and more usefully
    \[
    r_{k,D}(\sN) \geq \Expect{\vt,\va}{|A|}- \Expect{\vt,\va}{|T|},
    \]
with the expectation referring to choosing $\vt,\va$ uniformly from the torus $\TT^d$.
We have
    \[
    \Expect{\vt,\va}{|A|} = \Expect{\vt}{\Expect{\va}{|A|}} = \Expect{\vt}{N \vol(\annuli)}= N \vol(\annuli).
    \]
We also have
    \begin{equation*}
     \Expect{\vt,\va}{|T|}
        \leq \Expect{\vt,\va}{\sum E(D',a,b)} = \sum \Expect{\vt,\va}{E(D',a,b)}
     \end{equation*}
where $E(D',a,b)$ is 1 if $A$ contains a progression of type $(D',a,b)$, and is 0 otherwise, and the summation has $\type_{k,D}(\sN)$ summands. Using the assumption that $A_0$ is free of $k$-term $2D$-progressions, we are able to bound
    \[\Expect{\vt,\va}{E(D',a,b)}\]
efficiently in terms of the volume of $\annuli$ and the volume of a small sphere. We arrive at
     \begin{align*}
        \Expect{\vt,\va}{|T|} &\leq \type_{k,D}(\sN) \vol(\annuli)\vol(\ball),
     \end{align*}
which gives us a lower bound on $r_{k,D}(\sN)$ in terms of $\Psi,N_0,d,\delta$ and $A_0$. The work \cite{OBryant} gives a lower bound on the size of $A_0$, and optimization of the remaining parameters yields the result.

\section{Proof of Theorem~\ref{thm:main}}\label{sec:Proof}
The open interval $(a-b,a+b)$ of real numbers is denoted $a\pm b$. The interval $[1,N]\cap\ZZ$ of natural numbers is denoted $[N]$. The box $(\pm 2^{-D-1})^d$, which has Lebesgue measure $2^{-dD}$, is denoted $\BOX_D$. We define $\BOX_0=[-1/2,1/2)^d$.

Although we make no use of this until the very end of the argument, we set
    \[
    d := \left\lfloor 2^{n/2}\left( \frac{ \log \Psi(N) }{D}\right)^{1/(n+1)}  \right\rfloor.
    \]

Given $\vec x \in \RR^d$, we denote the unique element $\vec y$ of $\BOX_0$ with $\vec x - \vec y \in \ZZ^d$ as $\vec x \bmod \vec 1$.

A point $\vec x=\langle X_1,\dots, X_d \rangle$ chosen uniformly from $\BOX_D$ has components $X_i$ independent and uniformly distributed in $(-2^{-D-1},2^{-D-1})$. Therefore, $\|\vec x\|_2^2 = \sum_{i=1}^d X_i^2$ is the sum of $d$ iidrvs, and is consequently normally distributed as $d\to\infty$. Further, $\|\vec x \|_2^2$ has mean $\mu:=2^{-2D} d/12$ and variance $\sigma^2:=2^{-4D}d/180$.

Let $A_0$ be a subset of $[N_0]$ with cardinality $r_{k,2D}([N_0])$ that does not contain any $k$-term $2D$-progression, and assume $2\delta N_0\leq 2^{-2D}$. We define $\annuli$ in the following manner:
            \[\annuli :=\left\{\vec x \in \BOX_D \colon  \frac{\|\vec x\|_2^2-\mu}{\sigma} \in \bigcup_{a\in A_0} \left(z-\frac{a -1}{N_0} \pm \delta\right) \right\},\]
where $z\in \mu\pm\sigma$ is chosen to maximize the volume of $\annuli$. Geometrically, $\annuli$ is the union of $|A|$ spherical shells, intersected with $\BOX_D$. From \cite{OBryant}*{Lemma 3}, the Barry-Esseen central limit theorem and the pigeonhole principle yield:
\begin{lemma}[$\annuli$ has large volume] \label{lem:annuli}
            If $d$ is sufficiently large, $A_0\subseteq[N_0]$, and $2\delta\leq 1/n$, then the volume of $\annuli$ is at least
                \(\displaystyle \frac25 \, 2^{-dD} |A_0| \delta .\)
\end{lemma}

Set
    \[
    A:=A(\vt,\va) = \{n \in \sN \colon n \,\vt + \va \bmod\vec1 \in \annuli\},
    \]
which we will show is typically (with respect to $\vt,\va$ being chosen uniformly from $\BOX_0$) a set with many elements and few types of $D$-progressions. After removing one element from $A$ for each type of progression it contains, we will be left with a set that has large size and no $k$-term $D$-progressions.

Define $T:=T(\vt,\va)$ to be the set
\begin{equation*}
    \left\{a\in \sN \colon\;
    \begin{matrix}
        \text{$\exists b\in \RR, D'\in[D]$ such that $A(\vt,\va)$ contains} \\
        \text{a $k$-term progression of type $(D',a,b)$}
    \end{matrix}
        \right\},
\label{T.definition2}
\end{equation*}
which is contained in $A(\vt,\va)$. Observe that $A \setminus T$ is a subset of $\sN$ and contains no $k$-term $D$-progressions, and consequently
    \(r_{k,D}(\sN)\geq |A \setminus T|= |A| - |T|\)
for every $\vt,\va$. In particular,
    \begin{equation}\label{equ:rBound2}
    r_{k,D}(\sN) \geq \Expect{\vt,\va}{|A \setminus T|} =\Expect{\vt,\va}{|A| - |T|} = \Expect{\vt,\va}{|A|} - \Expect{\vt,\va}{|T|}.
    \end{equation}

First, we note that
    \begin{equation}\label{equ:ExpectA}
    \Expect{\vt,\va}{|A|} = \sum_{n\in\sN} \Prob{\vt,\va}{n\in A} = \sum_{n\in\sN} \Prob{\va}{n\in A}  = N \vol(\annuli).
    \end{equation}

Let $E{(D',a,b)}$ be 1 if $A$ contains a $k$-term progression of type $(D',a,b)$, and \mbox{$E{(D',a,b)}=0$} otherwise. We have
    \[|T| \leq \sum_{(D',a,b)} E{(D',a,b)},\]
where the sum extends over all types $(D',a,b)$ for which $D'\in[D]$ and there is a $k$-term $D'$-progression of that type contained in $\sN$; by definition there are $AP_{k,D}(\sN)$ such types.

Suppose that $A$ has a $k$-term progression of type $(D',a,b)$, with $D'\in[D]$. Let $p$ be a degree $D'$ polynomial with lead term $p_{D'}=b/D'!\not=0$, and $p(1),\dots,p(k)$ a $D'$-progression contained in $A$. Then
    \[\vec x_i:=p(i)\, \vt +\va \bmod\vec1 \in \annuli \subseteq \BOX_D.\]
We now pull a lemma from \cite{OBryant}*{Lemma 2}.
\begin{lemma}\label{lem:RL2} Suppose that $p(j)$ is a polynomial with degree $D'$, with $D'$-th coefficient $p_D'$, and set $\vec x_{j} := \vt\, p(j)+\va \mod \vec1$. If $\vec x_{1},\vec x_2,\dots, \vec x_k $ are in $\BOX_D$ and $k  \geq D+2$, then there is a vector polynomial $\vec P(j) = \sum_{i=0}^{D'} \vec P_i j^i$ with $\vec P(j)=\vec x_j$ for $j\in[k]$, and $D'!\vec P_{D'} =  \vt \,D'! p_{D'} \bmod \vec 1$.
\end{lemma}
Thus, the $\vec x_i$ are a $D'$-progression in $\RR^d$, say $\vec P(j)=\sum_{i=0}^{D'} \vec P_i j^i$ has $\vec P(j)=\vec x_j$ and $D'! \vec P_{D'}=D'! p_{D'} \,\vt \bmod\vec 1= b\, \vt \bmod \vec1$. Recalling that $z$ was chosen in the definition of $\annuli$, by elementary algebra
    \[Q(j):=\frac{ \|\vec P(j)\|_2^2 - \mu}{\sigma} - z\]
is a degree $2D'$ polynomial in $j$ (with real coefficients), and since $\vec P(j)=\vec x_j \in \annuli$ for $j\in [k]$, we know that
    \[Q(j) \in \bigcup_{a\in A_0} \left(-\frac{a-1}{N_0}\pm \delta\right)\]
for all $j\in [k]$, and also $Q(1),\dots,Q(k)$ is a $2D'$-progression. Define the real numbers $a_j \in A_0$, $\epsilon_j \in \pm \delta$ by
    \[Q(j) = - \frac{a_j-1}{N_0} + \epsilon_j.\]

For a finite sequence $(a_i)_{i=1}^k$, we define the forward difference $\Delta (a_i)$ to be the slightly shorter finite sequence $(a_{v+1}-a_v)_{v=1}^{k-1}$. The formula for repeated differencing is
            \[\Delta^m (a_i) =\left( \sum_{i=0}^{m} \binom{m}{i}(-1)^{i} a_{i+v}\right)_{v=1}^{k-m}.\]
We note that a nonconstant sequence $(a_i)$ with at least $2D+1$ terms is a $2D$-progression if and only if $\Delta^{2D+1}(a_i)$ is a sequence of zeros. If $a_i=p(i)$, with $p$ a polynomial with degree $2D$ and lead term $p_{2D}\not=0$, then $\Delta^{2D} (a_i) = ((2D)! p_{2D})$, a nonzero-constant sequence. Note also that $\Delta$ is a linear operator. Finally, we make use of the fact, provable by induction for $1\leq m \leq k$, that
            \[|\Delta^m (a_i)| \leq  2^{m-1} \left( \max_i a_i - \min_i a_i\right).\]

We need to handle two cases separately: either the sequence $(a_i)$ is constant or it is not. Suppose first that it is not constant. Since $a_i\in A_0$, a set without $k$-term $2D$-progressions, we know that $\Delta^{2D+1}(a_i)\not=(0)$, and since $(a_i)$ is a sequence of integers, for some $v$
    \[|\Delta^{2D+1}(a_i)(v)|\geq 1.\]
Consider:
    \[(0)=\Delta^{2D+1}(Q(i)) ={ \frac{1}{N_0} \Delta^{2D+1}(a_i) + \Delta^{2D+1}(\epsilon_i)}, \]
whence
    \[ |\Delta^{2D+1}(\epsilon_i)(v)| = \frac{1}{N_0} |\Delta^{2D+1}(a_i)(v)| \geq \frac{1}{N_0}.\]
Since $|\epsilon_i|<\delta$, we find that
    \[
    |\Delta^{2D+1}(\epsilon_i)(v)| = \left|\sum_{i=0}^{2D+1} \binom{2D+1}{i} (-1)^i \epsilon_{i+v} \right|
    < 2^{2D+1}\delta,
    \]
and since we assumed that $2\delta N_0\leq 2^{-2D}$, we arrive at the impossibility
    \[ \frac{1}{N_0} \leq |\Delta^{2D+1}(\epsilon_i)(v)| < 2^{2D+1} \delta \leq  2^{2D}\,\cdot\,\frac{2^{-2D}}{N_0} = \frac{1}{N_0}.\]

Now assume that $(a_i)$ is a constant sequence, say $a:=a_i$, so that
    \[Q(j) \in -\frac{a-1}{N_0} \pm \delta \]
for all $j\in[k]$. This translates to
    \[\|\vec P(j)\|_2^2 \in  \mu-(z-\frac{a-1}{N_0})\sigma \pm \delta\sigma.\]
Clearly a degree $2D'$ polynomial, such as $\|\vec P(j)\|_2^2$, cannot have the same value at $2D'+1$ different arguments; we pull now another lemma from \cite{OBryant}*{Lemma 1} that quantifies this.

\begin{lemma}\label{lem:RL1} Let $\delta,r$ be real numbers with $0\leq \delta \leq r$, and let $k,D$ be integers with $D\geq 1, k\geq 2D+1$. If $\vec P(j)$ is a polynomial with degree $D$, and $r-\delta \leq \| \vec P(j) \|_2^2 \leq r+\delta$ for $j\in[k]$, then the lead coefficient of $\vec P$ has norm at most $2^D \left.{(2D)!}\right.^{-1/2} \, \sqrt{\delta}$.
\end{lemma}

Using Lemma~\ref{lem:RL1}, the lead coefficient $\vec P_{D'}$ of $\vec P(j)$ satisfies
    \begin{equation*}
    \|D'! \vec P_{D'}\|_2 \leq D'!\,2^{D'} {(2D')!}^{-1/2} \sqrt{\delta \sigma} \leq \sqrt{ F \sigma \delta },
    \end{equation*}
where $F$ is an explicit constant. We have deduced that $E{(D',a,b)}=1$ only if
    \[
    a\,\vt+\va \bmod 1 \in \annuli \quad \text{and} \quad \|b\, \vt  \bmod 1\|_2 \leq \sqrt{ F \sigma \delta }.
    \]
Since $\va$ is chosen uniformly from $\BOX_0$, we notice that
    \[
    \Prob{\va}{a\,\vt+\va\bmod 1 \in \annuli}=\vol \annuli,
    \]
independent of $\vt$. Also, we notice that the event $\{\|b\, \vt \bmod 1\|_2 \leq \sqrt{ F \sigma \delta }\}$ is independent of $\va$, and that since $b$ is an integer, $\vt\bmod\vec 1$ and $b\,\vt\bmod \vec 1$ are identically distributed. Therefore, the event $\{\|b\, \vt \bmod 1\|_2 \leq \sqrt{ F \sigma \delta }\}$ has probability at most
    \[\vol \ball(\sqrt{ F \sigma \delta }) = \frac{2\pi^{d/2} (\sqrt{ F \sigma \delta })^{d}}{\Gamma(d/2) d},\]
where $\ball(x)$ is the $d$-dimensional ball in $\RR^d$ with radius $x$. It follows that
    \begin{equation*}
    \Prob{\vt,\va}{E{(D',a,b)}=1}\leq \vol \annuli \cdot\vol \ball(\sqrt{ F \sigma \delta }),
    \end{equation*}
and so
    \begin{equation}\label{equ:ExpectT}
    \Expect{\vt,\va}{|T|} \leq \type_{k,D}(\sN) \vol \annuli  \cdot \vol \ball(\sqrt{ F \sigma \delta }).
    \end{equation}

Equations~\eqref{equ:rBound2}, \eqref{equ:ExpectA}, and \eqref{equ:ExpectT} now give us
    \[\frac{r_{k,D}(N)}{N} \geq \vol(\annuli) \left( 1- \frac{\type_{k,D}(\sN)}{N} \vol \ball(\sqrt{ F \sigma \delta })\right).\]
Setting
    \[
    \delta = \frac{2e d }{\pi  F \sigma } \left( \frac{d}{d+2}\right)^{2/d} \frac{\Gamma(d/2)^{2/d}}{2ed} \left(\frac{\type_{k,D}(\sN)}{N}\right)^{-2/d}
    \sim C \frac{d^{1/2}}{\Psi(N)^{2/d}}
    \]
we observe that
    \[
    1- \frac{\type_{k,D}(\sN)}{N} \vol \ball(\sqrt{ F \sigma \delta }) = \frac {d}{d+2} \sim 1.
    \]
Now,
    \begin{align*}
    \frac{r_{k,D}(\sN)}{N}
        & \geq \vol \annuli \, \frac{d}{d+2}\\
        &\gg  2^{- dD} \,\delta |A_0|\\
        &\gg 2^{-dD} d^{1/2} \Psi(N)^{-2/d} |A_0|\\
        &= C \exp\left( -dD -\frac 2d \log \Psi(N) + \frac 12 \log d +\log|A_0|\right).
    \end{align*}
Recall that we set
    \[
    d := \left\lfloor 2^{n/2}\left( \frac{ \log \Psi(N) }{D}\right)^{1/(n+1)}  \right\rfloor.
    \]
If $2D<k\leq 4D$, we take $N_0=1$ and $A_0=\{1\}$ to complete the proof. If $k>4 D$, we set
    \[
    N_0 := C \frac{\Psi(N)^{2/d}}{d^{1/2}},
    \]
and use the bound
    \[
    |A_0|=r_{k,2D}(N_0) \geq C N_0 \exp\left( -n 2^{(n-1)/2} (2D)^{(n-1)/n} (\log N_0)^{1/n} + \frac{1}{2n}\log\log N_0\right),
    \]
proved in \cite{OBryant}, to complete the proof.

\section{Unanswered questions}\label{sec:questions}

Kolountzakis [personal communication] asks whether
    \[
    r_{3,1}([N]) = \min\{r_{3,1}(\sN) \colon \sN\subseteq \ZZ, |\sN|=N\}.
    \]
More generally, which set $\sN$ (for fixed $k,D,N$) minimizes $r_{k,D}(\sN)$? It is not even clear to this author which set maximizes $\type_{k,D}(\sN)$, nor even what that maximum is, although the interval $[N]$ is the natural suspect and has $\type_{k,D}([N])\leq 2^{D+1}N^2$.

We doubt that there is a subset of the squares with positive relative density that does not contain any 3-term arithmetic progressions, but haven't been able to prove such. We note that there are 4-term 2-progressions of positive cubes: $3^3, 16^3, 22^3, 27^3$ is the image of $0,1,2,3$ under $Q(x)=\frac{2483}{2}x^2+\frac{5655 }{2}x+27$. For which $k,D,p$ are there $k$-term $D$-progressions of perfect $p$-th powers, and when they exist how many types are there?

\end{document}